\documentclass{amsart}
\usepackage{amscd, amsfonts, amsmath, amsthm, amssymb,epsf}

\usepackage{tabularx}
\usepackage{ltablex}
\usepackage{longtable}
\usepackage{multirow}
\usepackage{hhline}
\usepackage{fancyhdr}

\usepackage{stmaryrd}
\usepackage[all]{xy}
\usepackage{srcltx}

\newcommand{\cal}{\mathcal}

\newcommand{\F}{\mathbb{F}}
\newcommand{\Q}{\mathbb{Q}}

\newcommand{\E}{\mathbb{E}}

\begin{document}

\author{Dinesh S. Thakur}

\title[Heuristics for expected  Mordell-Weil rank]{Expected  Mordell-Weil rank heuristics through Sato-Tate, Birch and Swinnerton-Dyer conjectures}

\address{ Department of Mathematics, University of Rochester, 
Rochester, NY 14627, dinesh.thakur@rochester.edu
}


\begin{abstract}
We present  an heuristic argument for the prediction of expected Mordell-Weil rank of elliptic curves over number fields, using Birch and Swinnerton-Dyer's 
original conjecture and Sato-Tate conjectures. We do calculations in some cases and raise questions about their relations, if any,  with the predictions of various average
rank models that have been considered. 
\end{abstract}


\maketitle



For an elliptic curve $E$ over the rational number field $\Q$,  write $N_p$ for the cardinality of $E(\F_p)$, and 
 let $P_x= \prod_{p\leq x}  N_p/p$.  

Then the strong original form of Birch and Swinnerton-Dyer (BSD) conjecture \cite{BSD} says that there is a positive constant $C=C(E)$ such that $P_x$ is asymptotic to 
$C(\log x)^r$, as $x$ tends to infinity, where $r$ is the (Mordell-Weil)  rank of $E$, i.e., the rank of the finitely generated group $E(\Q)$. 

Goldfeld \cite{G} showed that this implies the usual $L$-function form of this conjecture as well as the Riemann hypothesis for this function.
The finite field-function field Riemann Hypothesis, i.e.,  Hasse bounds imply much weaker result $\log(P_x) =O(x^{1/2+\epsilon})$, for any $\epsilon >0$, instead of the asymptotics $r\log\log(x)+B$ expected from above. But they  also imply that we can write $N_p=p+1-2\sqrt{p}\cos(\theta_p)$, $\theta_p\in [0, \pi]$ 

The distribution of  $\theta_p$ is predicted by the Sato-Tate (ST) conjecture \cite{Tate}  as equidistribution wrt. $(2/\pi) \sin^2(\theta)d\theta$,  for $E$ with no complex multiplications (non-CM).  This is now a theorem \cite{BGHT} (\cite{Taylor} for such $E$ with non-integral $j$-invariant.)

Now let us first calculate the expected value of $\log(P_x)$  under the Sato-Tate distribution on $\theta_p$ for non-CM elliptic curve. We write $\mathbb{E}$ for the expectation. 

Using $\log(1+y)=y-y^2/2+\cdots$, and   $\sum_{p\leq x} 1/p = \log\log(x)+B +O(1/\log(x))$, for some constant $B$,  we get

\begin{align*}
 \E(\log(P_x)) &= \sum \E(\log(1+1/p - 2 \cos(\theta_p)/\sqrt{p}))=\log\log(x)-\frac{1}{2}\log\log(x)+O(1)\\
  &=\frac{1}{2}\log\log(x)+O(1),
  \end{align*}
using   $\int_0^{\pi}2 \cos(\theta)\sin^2(\theta) d\theta =0$ and $\int_0^{\pi}2 \cos^2(\theta)\sin^2(\theta) d\theta =\pi/4$ to get the first two terms from the first 2 terms of the Taylor expansion, and putting  the error from the remaining convergent series in the big O term.

Comparison with the strong original BSD conjecture asymptotics of $r\log\log(x)$ implies heuristically  that {\it the expected value of the Mordell-Weil rank of  $E(\Q)$ in non-CM case is $1/2$}. 

\newpage

{\bf Remarks} (1) Note that we have not taken averages or expectations of actual ranks under any family, but just matched the asymptotics of the expectation of the quantity $\log(P_x)$ under ST for non-CM case with the asymptotics conjectured (a conjecture is, in fact, basically  just an `expectation'!) and called the result  heuristic expected value just because it comes from such a matching (and \cite{B} connection with the averages mentioned below).  Interestingly, this seems to fit well with the predictions (and some evidence and bounds in some cases) by Dorian Goldfeld \cite[Conjecture B]{G1} (for quadratic twists) and others (e.g., \cite[Conjecture 1.2]{PR}
over global fields) for the average rank as $1/2$ in various `averages scenarios'  when you take averages by ordering the elliptic curves by the heights (normalized coefficients sizes) or use families of quadratic twists ordered by the twist parameter.   (I do not know of any predictions (or rather results?) specifically mentioning  ordering by minimal discriminants or conductors, 
apart from Katz-Sarnak function field work\cite{KS}). We refer to \cite{BS, PPVW, P, PBour} for some results and surveys. Is this a coincidence or is there  a better 
explanation? In the Sato-Tate original case \cite{Tate}, one fixes $E$ and averages over primes, but as Birch \cite{B} proved, one gets the same probability distribution when we fix $p$ and average over curves $E$ over $\F_p$. So there might be an explanation of the coincidence. 

Andrew Sutherland suggested that probably the  ordering of $E$'s over $\Q$  by naive height has the property that for any
fixed $p$ and all sufficient large $M$ the first $M$ of the $E$'s  have
reductions mod $p$ that are equidistributed over the isomorphism classes
of $E$ over $\F_p$, and Birch might apply to connect with  actual averages. 

Peter Sarnak suggested that going from one curve (as in the expectation calculation) to a family, a serious complicated issue is how big the parameter $x$ is compared to the conductors 
of the family one is supposed to be averaging over.

 (2) Similar calculation gives the same expectation $1/2$ for non-CM elliptic curves over any number field,  by exactly the same argument and using that sum of the reciprocals of the norms of primes (up to $x$) still grows like $\log\log x$. Only the constants seem to depend on the number field. In contrast to \cite[Conjecture 1.2]{PR} for global fields, 
Dorian Goldfeld suggested to the author that $1/2$ may  not now be expected to be the average rank for the elliptic curves over number fields with many complex places, since  the discrepancy is known 
\cite{DD??}  for some quadratic twist families. The `expectations' and predictions by experts of the rank distributions have changed a few times over the time. 
The `Infinite family Rank  at most 21' heuristics of   Granville \cite{W} and \cite{PPVW} for elliptic curves over $\Q$ also seems to go wrong when generalized naively for general number fields, as pointed out in  \cite[Sec. 12]{PPVW}, \cite[Sec. 4]{P}.  

Next we look at CM  $E/Q$, a similar calculation now with the predicted uniform distribution of $\theta_p$ for ordinary primes and $a_p=0$ for super-singular p gives 

$\E(\log(P_x))=1/2 \log\log (x) +O(1)$,  now $1/2$ coming from the 
density of super-singular primes, the ordinary primes contribution being O(1) term: 
as $(1/2)\log\log(x)$ from the $y$ term (since $\int_0^{\pi}\cos(\theta) d\theta =0$) cancels from the $(1/2)\log\log(x)$ from $-y^2/2$ term
(since $\int_0^{\pi} 2 \cos^2(\theta)d\theta= \pi$.)  

Again,  this heuristically gives the expected rank as 1/2. 
This argument shows that the  expected rank is 1/2 for CM $E$ over field not containing the CM field and 0 for CM $E$ over field containing the CM field. 

 I do not know whether these are the `expected' averages in these cases. In fact, I could not find any published predictions. I would love to know references, if any. 

 Let us now calculate higher moments $\E((\log P_x)^n)$'s.

Write $x_i$ for $1/p_i-2\cos(\theta_p)/\sqrt{p_i}$, where $p_i$ is $i$-th prime. Then $\E(x_i)=1/p_i$ and $\E(x_i^2)=1/p_i$, so that $\sum \E(x_i)=\log\log (x) +O(1)$ and 
$\sum \E(x_i^2/2)=(1/2)\log\log (x) +O(1)$ by the calculations done above.  Since third and higher powers of $x_i$ can be ignored as they give convergent sums, we have, 
up to O(1), $\E((\log P_x)^n)=\E((\sum \log(1+x_i))^n)=\E(\sum \prod_{i=1}^n (x_i-x_i^2/2))$ 
since  terms with some repeated indices $i_1=i_2$ again contribute bounded quantities as they involve at least $3/2$ powers of primes in denominator. 

Let us do $n=2$ case first. Then up to O(1) again, with $i\neq j$ below, we have 

\begin{align*}
\E((\log(P_x))^2) &= \E(\sum \log(1+x_i)\log(1+x_j))
                           = \sum \E((x_i-x_i^2/2)(x_j-x_j^2/2)\\
                           & = \sum \E(x_ix_j-x_ix_j^2/2-x_i^2x_j/2+x_i^2x_j^2/4)\\
                           &= \sum [1/(p_ip_j)-1/(2p_ip_j)-1/(2p_jp_i)+1/(4p_ip_j)]\\
                           &= \sum 1/(4p_ip_j)
                           =(\log\log(x))^2/4
                           \end{align*}

For general $n$, we get, by the simplification above that 

$\E((\log P_x)^n)=\sum \E((\prod x_{i})(\prod (-x_{j}^2/2)=\sum \prod \E(x_i)\E(-x_j^2/2)$, where we take k of $x_i$ terms and $n-k$ of $-x_j^2/2$ terms, so we get 
$$(\sum_{k=0}^n(-1)^k{n\choose k}/2^k)(\log\log(x))^n=(1-1/2)^n(\log\log(x))^n=(\log\log(x))^n/2^n$$

{\bf Remarks}   The expectations (in the generic non-CM case that we only worked out so far) now seem  `unexpected': $n$-th moment is basically $(\log\log x)/2^n$, rather than $(\log\log x)/2$, which would be `expected',  if it is supposed to be the just average of $r^n$, since  half  are expected to be of rank 0 and half of rank 1, in the widely believed conjecture. Are the expectations  not supposed  to be reasonable averages  or is there some other good explanation for higher moments discrepancy with average considerations? 

Finally, we mention another (unrelated to above) observation regarding the rank conjectures. In a beautiful work \cite{PH}, Bjorn Poonen showed (among other things) that Hilbert's tenth problem over ${\cal O}_K$ is undecidable for every number field $K$, if for every number field $K$, there exists elliptic curve $E$ over $\Q$ with MW rank of $E$ over $\Q$ and over $K$ being one. But the optimistic conjecture \cite[2.6]{PH} that the if part is always true is incompatible with the BSD conjecture, since for some number fields K, e.g., $K=\Q(\sqrt{-1}, \sqrt{17})$
it is  known (by Karl Rubin see e.g., \cite[Sec. 8.1]{D12?}) that  the analytic rank of $E$ over $K$ is always even for all $E$ over $\Q$, by root number calculations. Bjorn Poonen tells me that he and Alexandra Shlapentokh proved later that for this application to Hilbert's tenth problem, it is enough to know that for
every quadratic  extension of number fields $L/K$, there is an elliptic curve $E$ of a positive rank over $K$ (not necessarily 1) that has the same rank over $L$. 
In fact, Karl Rubin and Barry Mazur \cite{MR} have proved that such $E$ exists, assuming that the Tate-Shafarevich group for elliptic curves is always finite.

\vskip .2truein
{\bf Acknowledgments} I thank Dorian Goldfeld, Peter Sarnak, Bjorn Poonen (especially for pointing out a wrong comment), Andrew Sutherland for their comments.

\end{document}